\documentclass[11pt]{article}
\usepackage{amsmath, upgreek}
\usepackage{amssymb}
\usepackage{color}
\usepackage{amscd}
\usepackage{xspace}
\usepackage{verbatim}
\usepackage{graphicx}
\renewcommand{\theequation}{\thesection.\arabic{equation}
}

\title{  Global $W^{1,p}$ regularity  for elliptic problem with measure source
and Leray-Hardy potential}
\author{{  Huyuan Chen\footnote{\noindent Department of Mathematics, Jiangxi Normal University,
Nanchang 330022, China. E-mail: chenhuyuan@yeah.net}} \\[3mm]
 {  Hichem Hajaiej\footnote{\noindent
California State University, Los Angeles, 5151, USA.  E-mail: hichem.hajaiej@gmail.com}}
}

\date{}
\begin{document}
 \maketitle


\newcommand{\txt}[1]{\;\text{ #1 }\;}
\newcommand{\tbf}{\textbf}
\newcommand{\tit}{\textit}
\newcommand{\tsc}{\textsc}
\newcommand{\trm}{\textrm}
\newcommand{\mbf}{\mathbf}
\newcommand{\mrm}{\mathrm}
\newcommand{\bsym}{\boldsymbol}
\newcommand{\scs}{\scriptstyle}
\newcommand{\sss}{\scriptscriptstyle}
\newcommand{\txts}{\textstyle}
\newcommand{\dsps}{\displaystyle}
\newcommand{\fnz}{\footnotesize}
\newcommand{\scz}{\scriptsize}
\newcommand{\be}{\begin{equation}}
\newcommand{\bel}[1]{\begin{equation}\label{#1}}
\newcommand{\ee}{\end{equation}}
\newcommand{\eqnl}[2]{\begin{equation}\label{#1}{#2}\end{equation}}
\newcommand{\barr}{\begin{eqnarray}}
\newcommand{\earr}{\end{eqnarray}}
\newcommand{\bars}{\begin{eqnarray*}}
\newcommand{\ears}{\end{eqnarray*}}
\newcommand{\nnu}{\nonumber \\}
\newtheorem{subn}{\name}
\renewcommand{\thesubn}{}
\newcommand{\bsn}[1]{\def\name{#1}\begin{subn}}
\newcommand{\esn}{\end{subn}}
\newtheorem{sub}{\name}[section]
\newcommand{\dn}[1]{\def\name{#1}}   
\newcommand{\bs}{\begin{sub}}
\newcommand{\es}{\end{sub}}
\newcommand{\bsl}[1]{\begin{sub}\label{#1}}
\newcommand{\bth}[1]{\def\name{Theorem}
\begin{sub}\label{t:#1}}
\newcommand{\blemma}[1]{\def\name{Lemma}
\begin{sub}\label{l:#1}}
\newcommand{\bcor}[1]{\def\name{Corollary}
\begin{sub}\label{c:#1}}
\newcommand{\bdef}[1]{\def\name{Definition}
\begin{sub}\label{d:#1}}
\newcommand{\bprop}[1]{\def\name{Proposition}
\begin{sub}\label{p:#1}}

\newcommand{\aand}{\quad\mbox{and}\quad}
\newcommand{\M}{{\cal M}}
\newcommand{\A}{{\cal A}}
\newcommand{\B}{{\cal B}}
\newcommand{\I}{{\cal I}}
\newcommand{\J}{{\cal J}}
\newcommand{\D}{\displaystyle}
\newcommand{\RR}{ I\!\!R}
\newcommand{\C}{\mathbb{C}}
\newcommand{\R}{\mathbb{R}}
\newcommand{\Z}{\mathbb{Z}}
\newcommand{\N}{\mathbb{N}}
\newcommand{\T}{{\rm T}^n}
\newcommand{\cuad}{{\sqcap\kern-.68em\sqcup}}
\newcommand{\abs}[1]{\mid #1 \mid}
\newcommand{\norm}[1]{\|#1\|}
\newcommand{\equ}[1]{(\ref{#1})}
\newcommand\rn{\mathbb{R}^N}
\renewcommand{\theequation}{\thesection.\arabic{equation}}
\newtheorem{definition}{Definition}[section]
\newtheorem{theorem}{Theorem}[section]
\newtheorem{proposition}{Proposition}[section]
\newtheorem{example}{Example}[section]
\newtheorem{proof}{proof}[section]
\newtheorem{lemma}{Lemma}[section]
\newtheorem{corollary}{Corollary}[section]
\newtheorem{remark}{Remark}[section]
\newcommand{\bremark}{\begin{remark} \em}
\newcommand{\eremark}{\end{remark} }
\newtheorem{claim}{Claim}


\newcommand{\rth}[1]{Theorem~\ref{t:#1}}
\newcommand{\rlemma}[1]{Lemma~\ref{l:#1}}
\newcommand{\rcor}[1]{Corollary~\ref{c:#1}}
\newcommand{\rdef}[1]{Definition~\ref{d:#1}}
\newcommand{\rprop}[1]{Proposition~\ref{p:#1}}
\newcommand{\BA}{\begin{array}}
\newcommand{\EA}{\end{array}}
\newcommand{\BAN}{\renewcommand{\arraystretch}{1.2}
\setlength{\arraycolsep}{2pt}\begin{array}}
\newcommand{\BAV}[2]{\renewcommand{\arraystretch}{#1}
\setlength{\arraycolsep}{#2}\begin{array}}
\newcommand{\BSA}{\begin{subarray}}
\newcommand{\ESA}{\end{subarray}}
\newcommand{\BAL}{\begin{aligned}}
\newcommand{\EAL}{\end{aligned}}
\newcommand{\BALG}{\begin{alignat}}
\newcommand{\EALG}{\end{alignat}}
\newcommand{\BALGN}{\begin{alignat*}}
\newcommand{\EALGN}{\end{alignat*}}
\newcommand{\note}[1]{\textit{#1.}\hspace{2mm}}
\newcommand{\Proof}{\note{Proof}}
\newcommand{\qeda}{\hspace{10mm}\hfill $\square$}
\newcommand{\qed}{\\
${}$ \hfill $\square$}
\newcommand{\Remark}{\note{Remark}}
\newcommand{\modin}{$\,$\\[-4mm] \indent}
\newcommand{\forevery}{\quad \forall}
\newcommand{\set}[1]{\{#1\}}
\newcommand{\setdef}[2]{\{\,#1:\,#2\,\}}
\newcommand{\setm}[2]{\{\,#1\mid #2\,\}}
\newcommand{\mt}{\mapsto}
\newcommand{\lra}{\longrightarrow}
\newcommand{\lla}{\longleftarrow}
\newcommand{\llra}{\longleftrightarrow}
\newcommand{\Lra}{\Longrightarrow}
\newcommand{\Lla}{\Longleftarrow}
\newcommand{\Llra}{\Longleftrightarrow}
\newcommand{\warrow}{\rightharpoonup}
\newcommand{
\paran}[1]{\left (#1 \right )}
\newcommand{\sqbr}[1]{\left [#1 \right ]}
\newcommand{\curlybr}[1]{\left \{#1 \right \}}
\newcommand{
\paranb}[1]{\big (#1 \big )}
\newcommand{\lsqbrb}[1]{\big [#1 \big ]}
\newcommand{\lcurlybrb}[1]{\big \{#1 \big \}}
\newcommand{\absb}[1]{\big |#1\big |}
\newcommand{\normb}[1]{\big \|#1\big \|}
\newcommand{
\paranB}[1]{\Big (#1 \Big )}
\newcommand{\absB}[1]{\Big |#1\Big |}
\newcommand{\normB}[1]{\Big \|#1\Big \|}
\newcommand{\produal}[1]{\langle #1 \rangle}

\newcommand{\thkl}{\rule[-.5mm]{.3mm}{3mm}}
\newcommand{\thknorm}[1]{\thkl #1 \thkl\,}
\newcommand{\trinorm}[1]{|\!|\!| #1 |\!|\!|\,}
\newcommand{\bang}[1]{\langle #1 \rangle}
\def\angb<#1>{\langle #1 \rangle}
\newcommand{\vstrut}[1]{\rule{0mm}{#1}}
\newcommand{\rec}[1]{\frac{1}{#1}}
\newcommand{\opname}[1]{\mbox{\rm #1}\,}
\newcommand{\supp}{\opname{supp}}
\newcommand{\dist}{\opname{dist}}
\newcommand{\myfrac}[2]{{\displaystyle \frac{#1}{#2} }}
\newcommand{\myint}[2]{{\displaystyle \int_{#1}^{#2}}}
\newcommand{\mysum}[2]{{\displaystyle \sum_{#1}^{#2}}}
\newcommand {\dint}{{\displaystyle \myint\!\!\myint}}
\newcommand{\q}{\quad}
\newcommand{\qq}{\qquad}
\newcommand{\hsp}[1]{\hspace{#1mm}}
\newcommand{\vsp}[1]{\vspace{#1mm}}
\newcommand{\ity}{\infty}
\newcommand{\prt}{\partial}
\newcommand{\sms}{\setminus}
\newcommand{\ems}{\emptyset}
\newcommand{\ti}{\times}
\newcommand{\pr}{^\prime}
\newcommand{\ppr}{^{\prime\prime}}
\newcommand{\tl}{\tilde}
\newcommand{\sbs}{\subset}
\newcommand{\sbeq}{\subseteq}
\newcommand{\nind}{\noindent}
\newcommand{\ind}{\indent}
\newcommand{\ovl}{\overline}
\newcommand{\unl}{\underline}
\newcommand{\nin}{\not\in}
\newcommand{\pfrac}[2]{\genfrac{(}{)}{}{}{#1}{#2}}

\def\ga{\alpha}     \def\gb{\beta}       \def\gg{\gamma}
\def\gc{\chi}       \def\gd{\delta}      \def\ge{\epsilon}
\def\gth{\theta}                         \def\vge{\varepsilon}
\def\gf{\phi}       \def\vgf{\varphi}    \def\gh{\eta}
\def\gi{\iota}      \def\gk{\kappa}      \def\gl{\lambda}
\def\gm{\mu}        \def\gn{\nu}         \def\gp{\pi}
\def\vgp{\varpi}    \def\gr{\rho}        \def\vgr{\varrho}
\def\gs{\sigma}     \def\vgs{\varsigma}  \def\gt{\tau}
\def\gu{\upsilon}   \def\gv{\vartheta}   \def\gw{\omega}
\def\gx{\xi}        \def\gy{\psi}        \def\gz{\zeta}
\def\Gg{\Gamma}     \def\Gd{\Delta}      \def\Gf{\Phi}
\def\Gth{\Theta}
\def\Gl{\Lambda}    \def\Gs{\Sigma}      \def\Gp{\Pi}
\def\Gw{\Omega}     \def\Gx{\Xi}         \def\Gy{\Psi}

\def\CS{{\mathcal S}}   \def\CM{{\mathcal M}}   \def\CN{{\mathcal N}}
\def\CR{{\mathcal R}}   \def\CO{{\mathcal O}}   \def\CP{{\mathcal P}}
\def\CA{{\mathcal A}}   \def\CB{{\mathcal B}}   \def\CC{{\mathcal C}}
\def\CD{{\mathcal D}}   \def\CE{{\mathcal E}}   \def\CF{{\mathcal F}}
\def\CG{{\mathcal G}}   \def\CH{{\mathcal H}}   \def\CI{{\mathcal I}}
\def\CJ{{\mathcal J}}   \def\CK{{\mathcal K}}   \def\CL{{\mathcal L}}
\def\CT{{\mathcal T}}   \def\CU{{\mathcal U}}   \def\CV{{\mathcal V}}
\def\CZ{{\mathcal Z}}   \def\CX{{\mathcal X}}   \def\CY{{\mathcal Y}}
\def\CW{{\mathcal W}} \def\CQ{{\mathcal Q}}
\def\BBA {\mathbb A}   \def\BBb {\mathbb B}    \def\BBC {\mathbb C}
\def\BBD {\mathbb D}   \def\BBE {\mathbb E}    \def\BBF {\mathbb F}
\def\BBG {\mathbb G}   \def\BBH {\mathbb H}    \def\BBI {\mathbb I}
\def\BBJ {\mathbb J}   \def\BBK {\mathbb K}    \def\BBL {\mathbb L}
\def\BBM {\mathbb M}   \def\BBN {\mathbb N}    \def\BBO {\mathbb O}
\def\BBP {\mathbb P}   \def\BBR {\mathbb R}    \def\BBS {\mathbb S}
\def\BBT {\mathbb T}   \def\BBU {\mathbb U}    \def\BBV {\mathbb V}
\def\BBW {\mathbb W}   \def\BBX {\mathbb X}    \def\BBY {\mathbb Y}
\def\BBZ {\mathbb Z}

\def\GTA {\mathfrak A}   \def\GTB {\mathfrak B}    \def\GTC {\mathfrak C}
\def\GTD {\mathfrak D}   \def\GTE {\mathfrak E}    \def\GTF {\mathfrak F}
\def\GTG {\mathfrak G}   \def\GTH {\mathfrak H}    \def\GTI {\mathfrak I}
\def\GTJ {\mathfrak J}   \def\GTK {\mathfrak K}    \def\GTL {\mathfrak L}
\def\GTM {\mathfrak M}   \def\GTN {\mathfrak N}    \def\GTO {\mathfrak O}
\def\GTP {\mathfrak P}   \def\GTR {\mathfrak R}    \def\GTS {\mathfrak S}
\def\GTT {\mathfrak T}   \def\GTU {\mathfrak U}    \def\GTV {\mathfrak V}
\def\GTW {\mathfrak W}   \def\GTX {\mathfrak X}    \def\GTY {\mathfrak Y}
\def\GTZ {\mathfrak Z}   \def\GTQ {\mathfrak Q}

\font\Sym= msam10 
\def\SYM#1{\hbox{\Sym #1}}
\newcommand{\bdw}{\prt\Gw\xspace}
\date{}
\maketitle

\begin{abstract}
In this paper, we develop the Littman-Stampacchia-Weinberger duality approach to obtain global
$W^{1,p}$ estimates for a class of  elliptic problems involving Leray-Hardy operators and  measure  sources
in a distributional framework associated to a dual formulation with a specific weight function.

\end{abstract}

\noindent
  \noindent {\small {\bf Key Words}:   Leray-Hardy potential; Duality approach; Radon measure.  }\vspace{1mm}

\noindent {\small {\bf MSC2010}:  35J75,   35B99,  35D99. }
\vspace{1mm}
\hspace{.05in}
\medskip

\setcounter{equation}{0}
\section{Introduction}

Our concern in this paper  is to derive    global $W^{1,p}$ estimates for weak solutions of non-homogeneous Hardy problem
 \begin{equation} \label{eq 1.1}
 \arraycolsep=1pt\left\{
\begin{array}{lll}
 \displaystyle    \mathcal{L}_\mu u =  \nu \quad
   &{\rm in}\ \ {\Omega},\\[1.5mm]
 \phantom{   L_\mu \,   }
 \displaystyle  u= 0\quad  &{\rm   on}\ \, \partial{\Omega},
 \end{array}\right.
 \end{equation}
where   $N\geq 3$, $\Omega$ is a bounded $C^2$ domain containing the origin,
 $\mathcal{L}_\mu$
is the Leray-Hardy operator defined by $\mathcal{L}_\mu= -\Delta   +\frac{\mu}{|x|^{2 }}$,
$\mu\geq \mu_0:=-\frac{(N-2)^2}{4}$ and $\nu$ is a non-homogeneous  source.

Note that  the inverse square potential in  $\mathcal{L}_\mu$ makes this operator degenerate at the origin.
Due to this fact, challenging problems appear in the related Hardy semilinear elliptic equations.
The numerous underlying problems were a source of motivation for many researchers. In order to solve
(\ref{eq 1.1}) in the standard case, $\nu$ is a function,  many aspects have been developed : Hardy inequality, see the references \cite{BM,D,F,FM,P,OK},  variational methods  \cite{DGL,FF} for example,  the singularities of semilinear Hardy elliptic equation   by   \cite{BP,CC,C,Fa,G,GV,N,Wa}  and the references therein.  When $N\geq 3$, $\mu_0\leq \mu<0$,  Dupaigne in \cite{D} (also see \cite{BDT,BM}) studied the weak solutions for the equation $\mathcal{L}_\mu u=u^p+tf\ \ {\rm in}\ \, \Omega$,
  subject to the zero Dirichlet boundary condition, with $t>0,\, p>1$ and $f\in L^1(\Omega)$,  in the distributional sense that $u\in L^1(\Omega)$, $u^p\in L^1(\Omega, \rho dx)$  and
\begin{equation}\label{w 1.1}
 \int_\Omega u\mathcal{L}_\mu \xi \,dx=\int_\Omega u^p \xi \,dx+t\int_\Omega f \xi \,dx,\quad \forall \,  \xi\in C^{1,1}_0(\Omega),
\end{equation}
where $\rho(x)={\rm dist}(x,\partial\Omega)$.

However, the normal distributional identity (\ref{w 1.1}) does not seem to extend to address other types of isolated singular solutions, such as the solutions of
\begin{equation}\label{eq 1.01}
\mathcal{L}_\mu u= 0\quad {\rm in}\quad  \R^N\setminus \{0\}.
\end{equation}
 Indeed,  when $\mu\geq\mu_0$,  problem (\ref{eq 1.01})
has two  radially symmetric solutions with the explicit formulas given by
\begin{equation}\label{1.1}
 \Phi_\mu(x)=\left\{\arraycolsep=1pt
\begin{array}{lll}
 |x|^{\tau_-(\mu)}\quad
   &{\rm if}\quad \mu>\mu_0\\[1mm]
 \phantom{   }
-|x|^{\tau_-(\mu)}\ln|x| \quad  &{\rm   if}\quad \mu=\mu_0
 \end{array}
 \right.\qquad {\rm and}\qquad \Gamma_\mu(x)=|x|^{ \tau_+(\mu)},
\end{equation}
where
$$
 \tau_-(\mu)=-\frac{N-2}2-\sqrt{\mu-\mu_0}\quad{\rm and}\quad  \tau_+(\mu)=-\frac{ N-2}2+\sqrt{\mu-\mu_0},
$$
 $\tau_-(\mu)$ and $\tau_+(\mu)$ are the zero points of $\mu-\tau(\tau+N-2)=0$. To simplify the notations, we write $\tau_+$ and $\tau_-$ in what follows.
Obviously, $\Phi_\mu $ has a too  strong singularity at the origin to hope that $\Phi_\mu\in L^1_{loc}(\R^N)$ when $\mu\geq \frac{N^2+4}{2}$. 
Due to its singularity,  $\Phi_\mu$ can't be viewed as a fundamental solution expressed by Dirac mass at origin in the normal distributional identity  for $\mu\not=0$.

  To overcome this difficulty,  \cite{CQZ} introduces a new distributional identity  for   $\Phi_\mu $ as
  \begin{equation}\label{eq0030}
\int_{\R^N}\Phi_\mu \mathcal{L}_\mu^*\xi\, d\gamma_\mu =c_\mu\xi(0),\quad\forall\, \xi\in C_c^{1,1}(\R^N),
\end{equation}
where $\delta_0$ is the Dirac mass at the origin, here and in what follows, 
\begin{equation}\label{L}
d\gamma_\mu(x) =\Gamma_\mu(x) dx,\quad \mathcal{L}^*_\mu=-\Delta -\frac{2\tau_+ }{|x|^2}\,x\cdot\nabla,\quad
c_{\mu}=\left\{\arraycolsep=1pt
\begin{array}{lll}
 2\sqrt{\mu-\mu_0}\,|\BBS^{N-1}|\quad
   &{\rm if}\quad \mu>\mu_0,\\[1.5mm]
 \phantom{   }
|\BBS^{N-1}| \quad  &{\rm  if}\quad \mu=\mu_0
 \end{array}
 \right.
\end{equation}
 and $\BBS^{N-1}$ is the sphere of  the unit ball in $\R^N$.
 The distributional identity (\ref{eq0030}) could be seen as
 \begin{equation}\label{eq003}
 \mathcal{L}_\mu   \Phi_\mu  =c_{\mu}\delta_0 \quad   {\rm in}\quad\R^N.
\end{equation}
 From above distributional identity,  $\Phi_\mu $ is a fundamental solution of Hardy operator $\mathcal{L}_\mu$.

Inspired by the new dual identity (\ref{eq0030}),   semilinear Hardy problems involving Radon measure
\begin{equation} \label{eq 1.1-absorption}
 \arraycolsep=1pt\left\{
\begin{array}{lll}
 \displaystyle    \mathcal{L}_\mu u +g(u)=  \nu \quad
   &{\rm in}\ \ {\Omega},\\[1.5mm]
 \phantom{   L_\mu +g(u)\,   }
 \displaystyle  u= 0\quad  &{\rm   on}\ \, \partial{\Omega}
 \end{array}\right.
 \end{equation}
 have been considered in \cite{ChVe} (also \cite{ChVe1} when the origin is at boundary) in a measure framework for $\nu=\gn\lfloor_{\Gw^*}+k\gd_0$,
 where $\Omega^*=\Omega\setminus\{0\}$, $\gn\lfloor_{\Gw^*}\in \mathfrak{M}(\Omega^*; \Gamma_\mu)$, here {\it   $\mathfrak{M}(\Omega^*; \Gamma_\mu)$ is denoted by the set of Radon measures $\gn$  in $\Gw^*$ such that
 \bel{I-5}
 \|\gn\|_{\mathfrak{M}(\Omega^*; \Gamma_\mu)}=\myint{\Gw^*}{}\Gg_\gm d|\gn|:=\sup\left\{\myint{\Gw^*}{}\gz d|\gn|\, :\gz\in C_c(\Gw^*),\,0\leq\gz\leq\Gg_\gm\right\}<\infty.
\ee}
 We denote by $\overline{\mathfrak{M}}(\Omega; \Gamma_\mu)$ the set of measures which can be written under the form
\bel{I-7}
\nu=\gn\lfloor_{\Gw^*}+k\gd_0,
\ee
where $\gn\lfloor_{\Gw^*}\in  \mathfrak{M}(\Omega; \Gamma_\mu)$ and $k\in\R$.
 Let $\overline \Omega^*\!\!:=\overline\Gw\setminus\{0\}$ and
\begin{equation}\label{test}
\BBX_\gm(\Gw)=\left\{\xi\in C_0(\overline\Gw)\cap C^1(\overline\Gw^*):|x|\CL_\gm^*\xi\in L^\infty(\Gw)\right\}.
\end{equation}
We note that $C^{1,1}_0(\overline\Gw)\subset \BBX_\gm(\Gw)$.

\begin{definition}\label{weak definition}
We say that $u$ is a weak solution (or a very weak solution) of (\ref{eq 1.1-absorption}) with $\nu=\gn\lfloor_{\Gw^*}+k\gd_0\in \overline{\mathfrak{M}}(\Omega; \Gamma_\mu)$,  if $u\in
L^1(\Omega,\,|x|^{-1}d{\gg_\mu} )$,  $g(u)\in L^1(\Omega, d\gg_\mu)$  and
\begin{equation}\label{weak sense-}
\int_\Omega u \mathcal{L}^*_\mu\xi \,d{\gg_\mu}+\int_\Omega g(u)\xi\, d\gg_\mu=\myint{\Gw}{}\xi d(\Gg_\gm\nu)+c_\mu k\xi(0)\quad \text{for all }\;\xi\in \BBX_\gm(\Gw).
\end{equation}

Similarly, a function $u$ is called  a weak solution (or a very weak solution) of (\ref{eq 1.1}) with   $\nu=\gn\lfloor_{\Gw^*}+k\gd_0$,
if  $u$ is in
$L^1(\Omega,\,|x|^{-1}d{\gg_\mu} )$ and satisfies  (\ref{weak sense-}) with $g\equiv0$.
\end{definition}

Away from the origin, the operator $\mathcal{L}_\mu$ is uniformly elliptic, therefore  local regularity properties for second
elliptic equation could be applied directly locally  in $\Omega\setminus\{0\}$. The main goal of this paper  is
to study the global regularity of solutions of (\ref{eq 1.1}) and to apply the global regularity to
semilinear elliptic equation (\ref{eq 1.1-absorption}). Our main result  reads as follows.

\begin{theorem}\label{teo 1}
Assume that $N\geq 3$,  $\mu>\frac34\mu_0$,    $\nu\in \overline{\mathfrak{M}}(\Omega; \Gamma_\mu)$ with the form  $\nu=\gn\lfloor_{\Gw^*}+k\gd_0$
for $k\in \R$,  and $u$ is a  very weak solution of (\ref{eq 1.1}).
Let
\begin{equation}\label{critical 1}
 p^*_\mu=\frac{N+\tau_+}{N-1+\tau_+}.
\end{equation}
Then
$(i)$   $u\Gamma_\mu \in W^{1,p}_0(\Omega)$ with $p\in[1,\frac{N}{N-1})$ and there exists $C_1>0$ depending on $N,\ p,\, \mu$ and $\Omega$ such that
\begin{equation}\label{1.2-0}
\norm{u\Gamma_\mu  }_{W^{1,p}(\Omega)} \le C_1 (\norm{\nu}_{ \mathfrak{M}(\Omega^*; \Gamma_\mu)}+k).
\end{equation}

$(ii)$   $u  \in W^{1,p}_0(\Omega,d\gg_\mu)$ with $p\in[1, \min\{p^*_\mu, \frac{N}{N-1}\})$ and there exists $C_2>0$ depending on $N,\ p,\, \mu$ and $\Omega$ such that
\begin{equation}\label{1.2-0}
\norm{u   }_{W^{1,p}(\Omega, d\gg_\mu)} \le C_2 (\norm{\nu}_{ \mathfrak{M}(\Omega^*; \Gamma_\mu)}+k).
\end{equation}

\end{theorem}

\begin{remark} $(i)$ $p^*_\mu\leq \frac{N}{N-1}$ for $\mu\geq 0$ and  $p^*_\mu> \frac{N}{N-1}$
for $\mu\in[\mu_0, 0)$;

$(ii)$  An equivalent form of $p^*_\mu$ is $p^*_\mu=1-\frac{2}{\tau_-}$ by definition of $\tau_-$, and
$\min\{p^*_\mu, \frac{N}{N-1}\}$ is the critical exponent for the existence of (\ref{eq 1.1-absorption})
for $\nu\in\overline{\mathfrak{M}}(\Omega; \Gamma_\mu)$ in $\cite{ChVe}$, is the critical exponents for
isolated singular solutions of $\CL_\mu u=u^p$ in $\Omega\setminus\{0\}$, $u=0$ on $\partial\Omega$ in \cite{CZ}.

\end{remark}

 It is well-known that  the classical method for global  $W^{1,p}$ estimates for second order elliptic equations  named as  Littman-Stampacchia-Weinberger duality approach ( LSW   approach for short) in \cite{LSW,S},  with the coefficients being in $L^\infty$.  LSW approach has been used to deal with various elliptic equations, see  references \cite{BW,BP,CP,Fa,P}.    Global  $W^{1,p}$ estimates in Theorem \ref{teo 1} are derived by LSW approach.  Due to the inverse-square potential, the Hardy operator  is degenerated. Our dual operator  for  $\mathcal{L}_\mu$ is no longer $\mathcal{L}_\mu$, but is $\CL^*_\mu$ in our distributional sense, see Definition \ref{weak definition}. Observe that the coefficient of  dual operator $\mathcal{L}_\mu^* $ in gradient term is singular.
That is the main reason for assumptions that $N\geq3$ and $\mu>\frac{3}{4}\mu_0$, and under these assumptions,  
the solutions for nonhomogeneous Dirichlet problem
\begin{equation}\label{dual problem}
\mathcal{L}_\mu^* u=f_n+{\rm div}F_n\ \ {\rm in}\ \, \Omega,\qquad u=0\ \ {\rm on}\ \, \partial\Omega
\end{equation}
keep bounded 
 in the approximation of $f$ and $F$ in $L^r$ nonhomogeneous terms with $r>N$ by smooth functions $\{f_n\}_n$ and $\{F_n\}_n$. These uniform estimates
 enable us to make use of the LSW duality approach.

While if $\mu>0$ for $N=2$ or if $\mu\in [\mu_0,\frac34\mu_0]$ for $N\geq 3$,  dual problem (\ref{dual problem}) has a unique classical solution by constructing super and sub solutions
when $f$ and $F$ are regular. However, it is very challenging to obtain weak solutions and  $L^\infty$ estimates in the approximation of $f$ and  $F$ in $L^2$ spaces.
It seems that one has to come up with a new and self-contained
approach to solve this problem.

For our semilinear elliptic problem (\ref{eq 1.1-absorption}),  a direct corollary for a global regularity  could be obtained as follows.
 \begin{corollary}\label{cr 1}
 Assume that $N\geq 3$,  $\mu>\frac34\mu_0$ and  $\nu\in \overline{\mathfrak{M}}(\Omega; \Gamma_\mu)$.
 If problem  (\ref{eq 1.1-absorption}) admits a weak solution  $u_g$, then $u_g\Gamma_\mu\in W^{1,p}_0(\Omega)$ for $p\in[1, \frac{N}{N-1})$ and $u\in W^{1,p}_0(\Omega,d\gg_\mu)$ with $p\in[1, \min\{p^*_\mu, \frac{N}{N-1}\})$.
\end{corollary}

The global regularity  $u_g\Gamma_\mu\in W^{1,p}_0(\Omega)$  obtained in Corollary \ref{cr 1} indicates that
the zero Dirichlet  boundary restriction of   (\ref{eq 1.1-absorption}) could be understood in the trace sense.

\medskip

 The rest of the paper is organized as follows. In Section 2, we approximate the solution of
 $ \mathcal{L}_\mu^* v = f \Gamma_\mu  +  {\rm div} (F\Gamma_\mu )$,  subject to zero Dirichlet boundary condition and
  build $L^\infty$ estimate for functions $f,F$ in suitable spaces. Section  3 is devoted to prove Theorem \ref{teo 1} by developing LSW duality approach.

\setcounter{equation}{0}
\section{ Weak solution of Dual Hardy problem}

\subsection{Existence}

 The essential point to  obtain $W^{1,p}$ estimates by applying the Littman-Stampacchia-Weinberger duality method
 is to get  $L^\infty$ estimates of the related dual  nonhomogeneous problem.  In what follows,  we denote by $c_i$ a generic  positive constant.

 Let us consider  the  solution $w$ of the problem
\begin{equation}\label{eq 2.1}
\arraycolsep=1pt\left\{
\begin{array}{lll}
 \displaystyle   \mathcal{L}_\mu^* w = f   +  {\rm div}F\qquad
   &{\rm in}\quad  {\Omega},\\[2mm]
 \phantom{  L_\mu   }
 \displaystyle  w= 0\qquad  &{\rm   on}\quad \partial{\Omega},
 \end{array}\right.
\end{equation}
where $f:\Omega\to\R$ and $F:\Omega\to\R^N$.
Let  $\mathbb{H}^1_{0}(\Omega)$ be the Hilbert space defined as the closure of  functions in $C^\infty_c(\Omega)$
 under the norm
 $$\norm{u}_{\mathbb{H}^1_{0}(\Omega)}=\sqrt{<u,u>_{\mathbb{H}^1_{0}(\Omega)}}$$
with the inner product
$$<u,v>_{\mathbb{H}^1_{0}(\Omega)}=\int_\Omega\nabla u\cdot \nabla v\, dx.$$

{\it A function $w$ is called  a weak solution of (\ref{eq 2.1}) if $w\in \mathbb{H}^1_{0}(\Omega)$
 and
$$
 \int_\Omega   \Big(\nabla w\cdot \nabla v-  \frac{2\tau_+  }{ |x|^2 }\, (x\cdot \nabla w) v \Big)\, dx   =\int_{\Omega}(fv - F\cdot\nabla v)\, dx\quad \text{ for any } v\in \mathbb{H}^1_{0}(\Omega).
$$

A function $w$ is called  a classical solution of (\ref{eq 2.1}) if $w\in \mathbb{H}^1_{0}(\Omega)\cap C^2(\Omega\setminus\{0\})$
 and $w$ verifies equation (\ref{eq 2.1}) pointwisely in $\Omega\setminus\{0\}$.
}

Our  existence results state as follows.
 \begin{proposition}\label{pr 2.0}
Assume that $N\geq 3$ and $
 \mu>\frac34 \mu_0.$ Then\smallskip

 $(i)$ for $\theta\in(0,1)$, for any   $f \in C^\theta(\bar\Omega)$ and $F  \in  C^{1,\theta}(\bar\Omega;\R^N)$,
problem (\ref{eq 2.1}) has a bounded,  unique classical solution.

 $(ii)$  for $f \in L^2(\Omega)$ and $F  \in  L^2(\Omega;\R^N)$,
  problem (\ref{eq 2.1}) has   a unique weak solution $w\in \mathbb{H}^1_{0}(\Omega)$ satisfying that
 \begin{equation}\label{2.1-a1}
   \norm{w}_{\mathbb{H}^1_{0}(\Omega)}\le c_1\left(\norm{f }_{L^2(\Omega )} +\norm{ F }_{L^2(\Omega )} \right).
 \end{equation}
\end{proposition}

 Due to the singularity of the coefficient of $\CL^*_\mu$, we use the following sequence of operator
 to approximate it.
Given $\epsilon\in(0,\, 1)$, we denote
$$\mathcal{L}^*_{\mu,\epsilon}=-\Delta -2\tau_+ \frac{ 1_{B_\ge^c}}{  |x|^2  }\,x\cdot\nabla$$
and we consider the approximated problem
\begin{equation}\label{eq 2.1-ep}
\arraycolsep=1pt\left\{
\begin{array}{lll}
 \displaystyle   \mathcal{L}_{\mu,\epsilon}^* w = f   +  {\rm div}F\qquad
   &{\rm in}\quad  {\Omega},\\[2mm]
 \phantom{  L_\mu   }
 \displaystyle  w= 0\qquad  &{\rm   on}\quad \partial{\Omega},
 \end{array}\right.
\end{equation}
where $f:\Omega\to\R$, $F:\Omega\to\R^N$ and $1_{B_\ge^c}$ is the characteristic function in $\R^N\setminus B_\ge(0)$.

Since $2\tau_+  1_{B_\ge^c}   |x|^{-2}\,x$ is bounded,  smooth in $\Omega\setminus \partial B_\ge(0)$  and
$\mathcal{L}^*_{\mu,\epsilon}$ is a uniformly elliptic operator, then for any $f\in C^1(\bar \Omega)$
and $F\in C^2(\bar \Omega, \R^N)$, problem (\ref{eq 2.1-ep}) admits a unique classical solution $u_\ge $ by the
 Perron's super and sub solutions' method and comparison principle.
On the other hand, by the fact that $0\in \Omega$, there exist $\ge_0\in(0.\,1)$ and $R_0>\ge_0$ such that
$$B_{\ge_0}(0)\subset \Omega\subset B_{R_0}(0).$$

\begin{lemma}\label{lm 2.1-ep0}
Assume that $N\geq 2$, $\ge\in[0,\,\ge_0)$, $\mu>\frac34\mu_0$,  $f\in C^1(\bar \Omega)$,
 $F\in C^2(\bar \Omega, \R^N)$ and $u_\ge$ is the unique classical solution of problem (\ref{eq 2.1-ep}).
Then  there exists $c_2>0$ independent of $f, F$ and $\ge$ such that
 \begin{equation}\label{2.1-b0}
   \norm{\nabla u_\ge}_{L^2(\Omega)}\le c_2\left(\norm{f }_{L^2(\Omega )} +\norm{ F }_{L^2(\Omega )} \right).
 \end{equation}

\end{lemma}
{\bf Proof.} It is obvious that $u_\ge\in \mathbb{H}^1_{0}(\Omega)$ for $\ge\in(0,\, 1)$
and we observe that
 \begin{equation}\label{2.1-b2}
  (N-2)\tau_+>-\frac{(N-2)^2}{4}\ \text{ if } \ \mu>-\frac34\frac{(N-2)^2}{4},
  \end{equation}
  that is,
  $$\frac{4\tau_+ }{N-2}>-1.$$

Multiplying $u_\ge$ in  (\ref{eq 2.1-ep}) and integrating over $\Omega$, we obtain that
$$
 \int_\Omega   \Big(|\nabla u_\ge|^2-  \frac{\tau_+  1_{B_\ge^c} }{ |x|^2 }\, x\cdot \nabla (u_\ge^2)\Big)\, dx   =\int_{\Omega}(fu_\ge - F\cdot\nabla u_\ge)\, dx.
$$
By H\"{o}lder inequality and Poincar\'e inequality, we have that
\begin{eqnarray*}
\int_{\Omega}(fu_\ge  - F\cdot\nabla u_\ge)\, dx&\le&  (\norm{f}_{L^2(\Omega)}+\norm{F}_{L^2(\Omega;\R^N)}) (\norm{u_\ge}_{L^2(\Omega)}+\norm{\nabla u_\ge}_{L^2(\Omega)})
\\&\le& c_3 (\norm{f}_{L^2(\Omega)}+\norm{F}_{L^2(\Omega;\R^N)})  \norm{\nabla u_\ge}_{L^2(\Omega)},
 \end{eqnarray*}
   where $c_3>0$ is independent of $f, F$ and $\ge$.
   Observe that $\frac{x}{  |x|^2} u_\ge^2\in W^{1,1}(\Omega\setminus B_\ge(0):\R^N)$ and the divergence theorem, (see e.g. \cite[Theorem 6.3.4]{W}),  imply that
  \begin{eqnarray}
 \int_{\Omega\setminus B_\ge(0)}  \nabla (u_\ge^2) \frac{x}{ |x|^2} dx +\int_{\Omega\setminus B_\ge(0)}  u_\ge^2 {\rm div}(\frac{x}{  |x|^2} u)dx &=&  \int_{\Omega\setminus B_\ge(0)}  {\rm div}(\frac{x}{  |x|^2} u_\ge^2)dx \nonumber
\\&=&    \int_{\partial B_\ge(0)}  u_\ge^2 \frac{x}{  |x|^2}\cdot \frac{-x}{|x|} d\omega(x)\leq 0,\label{div 1}
 \end{eqnarray}
thus,
 \begin{eqnarray*}
 \int_\Omega   \Big(|\nabla u_\ge|^2-\frac{\tau_+  1_{B_\ge^c} }{ |x|^2 }\, x\cdot \nabla (u_\ge^2)\Big)\, dx &\geq &  \int_\Omega    |\nabla u_\ge|^2dx +\tau_+   \int_\Omega   u_\ge^2\, {\rm div}(\frac{x}{  |x|^2}1_{B_\ge^c}) dx
\\&=&  \norm{\nabla u_\ge}_{L^2(\Omega)}^2 +  (N-2)\tau_+\int_{\Omega\setminus B_\ge(0)}  \frac{u_\ge^2}{  |x|^2}dx.
 \end{eqnarray*}
When $\tau_+\geq 0$ i.e. $\mu\geq0$,
 \begin{eqnarray*}
 \int_\Omega   \Big(|\nabla u_\ge|^2-\frac{\tau_+  1_{B_\ge^c} }{ |x|^2 }\, x\cdot \nabla (u_\ge^2)\Big)\, dx  &\geq& \norm{\nabla u_\ge}_{L^2(\Omega)}^2.
 \end{eqnarray*}
 When $N\geq 3$ and $0>\mu>-\frac34\frac{(N-2)^2}{4}$, by Hardy inequality, we have that
 \begin{eqnarray*}
 \int_\Omega   \Big(|\nabla u_\ge|^2-  \frac{\tau_+  1_{B_\ge^c} }{ |x|^2 }\, x\cdot \nabla (u_\ge^2)\Big)\, dx     &\geq&  \norm{\nabla u_\ge}_{L^2(\Omega)}^2 + (N-2)\tau_+ \int_{\Omega}  \frac{u_\ge^2}{  |x|^2}dx
 \\&\geq&(1+\frac{4\tau_+ }{N-2})\norm{\nabla u_\ge}_{L^2(\Omega)}^2 .
 \end{eqnarray*}

 Therefore, there exists $c_4>0$ independent of $\ge, f$ and $F$ such that
 $$\norm{\nabla u_\ge}_{L^2(\Omega)}\leq c_4 (\norm{f}_{L^2(\Omega)}+\norm{F}_{L^2(\Omega;\R^N)}),$$
 which completes the proof.\hfill$\Box$ \medskip

\begin{lemma}\label{lm 2.1-ep1}
Assume that $N\geq 3$, $0<\ge_1<\ge_2<\ge_0$, $\mu>\frac34\mu_0$ and  $f\in C^1(\bar \Omega)$ and $F\in C^2(\bar \Omega, \R^N)$ and $u_{\ge_1},\, u_{\ge_2}$ are the unique classical solutions of problem (\ref{eq 2.1-ep}) with $\ge=\ge_1$ and $\ge=\ge_2$ respecitvely.

Then  there exists $c_5>0$ independent of   $\ge_1$ and $\ge_2$ such that
 \begin{equation}\label{2.1-b2}
   \norm{\nabla (u_{\ge_1}-u_{\ge_2})}_{L^2(\Omega)}\le c_5 \left(\norm{f }_{L^2(\Omega )} +\norm{ F }_{L^2(\Omega )} \right)\sqrt{\ge_2^{N-2}-\ge_1^{N-2}}.
 \end{equation}

\end{lemma}
\noindent{\bf Proof.}
Due to the boundedness of $\Omega$, there exists $R_0\geq \ge_0$ such that
$\Omega\subset B_{R_0}(0)$. Let $a_0=\|f\|_{L^\infty(\Omega)} +\|{\rm div} F\|_{L^\infty(\Omega)}$ and
\begin{equation}\label{2.1-b1}
W_0(x)=\frac{a_0}{N+2}(R_0^2-|x|^2),\quad\forall\, x\in\Omega.
\end{equation}
For $\mu>\frac34\mu_0$, we have that  $2N+2\tau_+>N+2$ and
$$\mathcal{L}^*_{\mu,\epsilon} W_0=\frac{a_0}{N+2}(2N+2\tau_+1_{B_\ge^c})\geq a_0. $$
In addition, Comparison Principle implies that for any $\ge\in(0,\, \ge_0)$,
$$|u_{\ge}|\leq W_0\quad {\rm in}\ \ \Omega.$$
Direct computation shows that
 \begin{eqnarray*}
 0&=&\int_\Omega  (\mathcal{L}_{\mu,\epsilon_1}^*u_{\ge_1} -\mathcal{L}_{\mu,\epsilon}^*u_{\ge_2} )(u_{\ge_1}-u_{\ge_2})\, dx
 \\&\geq &  \int_\Omega    |\nabla (u_{\ge_1}-u_{\ge_2})|^2dx +\tau_+   \int_{\Omega}   (u_{\ge_1}-u_{\ge_2})^2\, {\rm div}(\frac{x}{  |x|^2}1_{B_{\ge_1}^c}) dx
 \\&&-2\tau_+ \int_{  B_{\ge_2}(0)\setminus B_{\ge_1}(0)} (u_{\ge_1}-u_{\ge_2}) \frac{x\cdot \nabla u_{\ge_2}}{|x|^2}dx
\\&\geq &  (1+\frac{4\tau_+ }{N-2})\norm{\nabla  (u_{\ge_1}-u_{\ge_2}) }_{L^2(\Omega)}^2  - 2 |\tau_+| (\int_{B_{\ge_2}(0)\setminus B_{\ge_1}(0)} \frac{|u_{\ge_1}-u_{\ge_2}|^2 }{|x|^2}dx)^{\frac12} \norm{\nabla u_{\ge_2}}_{L^2(\Omega)}
\\&\geq &  (1+\frac{4\tau_+ }{N-2})  \norm{\nabla  (u_{\ge_1}-u_{\ge_2}) }_{L^2(\Omega)}^2  -4a_0 |\tau_+| (\int_{B_{\ge_2}(0)\setminus B_{\ge_1}(0)} \frac{1 }{|x|^2}dx)^{\frac12}\norm{\nabla u_{\ge_2}}_{L^2(\Omega)}   ,
 \end{eqnarray*}
 thus, by Lemma \ref{lm 2.1-ep0}, we have that
 $$ \norm{\nabla  (u_{\ge_1}-u_{\ge_2}) }_{L^2(\Omega)}^2 \leq  c_5 \left(\norm{f }_{L^2(\Omega)} +\norm{ F }_{L^2(\Omega)} \right)\sqrt{\ge_2^{N-2}-\ge_1^{N-2}},$$
 which completes the proof.\hfill$\Box$ \medskip

\begin{remark}
In the above Lemma, the dimension $N\geq 3$ plays an essential role in the calculation of $\int_{B_{\ge_2}(0)\setminus B_{\ge_1}(0)} \frac{1 }{|x|^2}dx$.

\end{remark}

\noindent{\bf Proof of Proposition \ref{pr 2.0}.}
$(i)$  Take $\ge_n=\frac1n$ with $n\in\N$ and $n\geq \frac1{\ge_0}$ and let $u_n$ be  the unique solution of
 (\ref{eq 2.1-ep}), then by Lemma \ref{lm 2.1-ep1}, we have that $\{u_n\}$ is a Cauchy sequence in $\mathbb{H}^1_{0}(\Omega)$.
  Therefore, there exists $u_0\in \mathbb{H}^1_{0}(\Omega)$ such that
 $$u_n\to u_0\quad{\rm in}\ \mathbb{H}^1_{0}(\Omega)\quad{\rm as}\ \, n\to+\infty.$$
 Multiplying $v\in \mathbb{H}^1_{0}(\Omega)$ in  (\ref{eq 2.1-ep}) and integrating over $\Omega$, we obtain that
 \begin{equation}\label{2.1-b2}
 \int_\Omega   \Big(\nabla u_n\cdot \nabla v-  \frac{2\tau_+  }{ |x|^2 }\, (x\cdot \nabla u_n) v 1_{B_{\frac1n}^c}\Big)\, dx   =\int_{\Omega}(fv - F\cdot\nabla v)\, dx.
 \end{equation}
 Passing to the limit  in the identity (\ref{2.1-b2}) as $n\to+\infty$, we have that
 \begin{equation}\label{2.1-b3}
 \int_\Omega   \Big(\nabla u_0\cdot \nabla v-  \frac{2\tau_+  }{ |x|^2 }\, (x\cdot \nabla u_0) v \Big)\, dx   =\int_{\Omega}(fv - F\cdot\nabla v)\, dx\quad \text{ for any } v\in \mathbb{H}^1_{0}(\Omega).
 \end{equation}
 From the proof of Lemma \ref{lm 2.1-ep1},  the sequence $\{u_n\}$ is uniformly bounded by barrier $W_0$ defined by (\ref{2.1-b1}), so $|u_0|\leq W_0$ a.e. in $\Omega$. Standard interior regularity shows that $u_0$ is a classical solution of (\ref{eq 2.1}).\smallskip

(ii)  For $f \in L^2(\Omega)$ and $F  \in  L^2(\Omega;\R^N)$,  Let  $\{f_m\}$ and $\{F_m\}$ be two sequence in  $C^1(\bar\Omega)$ and $C^2(\bar\Omega;\R^N)$ converging to $f$ and $F$ in $L^2(\R^N)$ and
$L^2(\Omega;\R^N)$ respectively.

By Lemma \ref{lm 2.1-ep0} with $\ge=0$,
\begin{equation}\label{2.1-b00}
   \norm{\nabla u_{m} }_{L^2(\Omega)}\le c_5\left(\norm{f_{m}  }_{L^2(\Omega )} +\norm{ F_{m}  }_{L^2(\Omega )} \right)
 \end{equation}
 and
 \begin{equation}\label{2.1-b0}
   \norm{\nabla (u_{m_1}-u_{m_2})}_{L^2(\Omega)}\le c_5\left(\norm{f_{m_1}-f_{m_2} }_{L^2(\Omega )} +\norm{ F_{m_1}-F_{m_2} }_{L^2(\Omega )} \right).
 \end{equation}
Therefore, we have that $\{u_m\}$ is a Cauchy sequence in $\mathbb{H}^1_{0}(\Omega)$ and there exists $u\in \mathbb{H}^1_{0}(\Omega)$ such that $u_m\to u$ in $\mathbb{H}^1_{0}(\Omega)$ as $m\to+\infty$.
Additionally, we have that
\begin{equation}\label{2.1-b4}
 \int_\Omega   \Big(\nabla u_m\cdot \nabla v-  \frac{2\tau_+  }{ |x|^2 }\, (x\cdot \nabla u_m) v \Big)\, dx   =\int_{\Omega}(fv - F\cdot\nabla v)\, dx\quad \text{ for any } v\in \mathbb{H}^1_{0}(\Omega).
 \end{equation}
Passing to the limit as $m\to+\infty$,  we have that
\begin{equation}\label{2.1-b5}
 \int_\Omega   \Big(\nabla u\cdot \nabla v-  \frac{2\tau_+  }{ |x|^2 }\, (x\cdot \nabla u) v \Big)\, dx   =\int_{\Omega}(fv - F\cdot\nabla v)\, dx\quad \text{ for any } v\in \mathbb{H}^1_{0}(\Omega).
 \end{equation}
 Therefore, $u$ is a weak solution of (\ref{eq 2.1}) and (\ref{2.1-a1}) follows by (\ref{2.1-b00}) directly for
 $\mu>-\frac34\frac{(N-2)^2}{4}$. \smallskip

 The uniqueness follows by Lemma \ref{lm 2.1-ep0}. In fact, if there exists two (classical or weak) solutions $u_1, u_2\in \mathbb{H}^1_{0}(\Omega)$ of (\ref{eq 2.1}), then it follows by Lemma \ref{lm 2.1-ep0} that
 $\norm{\nabla (u_1-u_2)}_{L^2(\Omega)}\leq 0$ and then $u_1=u_2$ a.e. in $\Omega$.\hfill$\Box$ \medskip

\subsection{$L^\infty$ esitmates for dual problems}

Given $k>0$, let $S_k$ be the function defined for $t\in \R$ by
\begin{equation}\label{2.0}
 S_k(t)=\left\{\arraycolsep=1pt
\begin{array}{lll}
 t+k\qquad
   &{\rm if}\quad \ t<-k,\\[1mm]
 \phantom{   }
0\qquad
   &{\rm if}\quad -k\le t\le k,\\[1mm]
 \phantom{   }
 t-k \qquad  &{\rm  if}\quad \ t>k.
 \end{array}
 \right.
\end{equation}

\begin{lemma}\label{lm 2.0}
Assume that $w$ is a measurable function in $L^1(\Omega)$ and there exist $\alpha>1$ and $A>0$ such that for every $k>0$,
\begin{equation}\label{2.3}
 \norm{  S_k(w) }_{L^1(\Omega)} \le A \{|w|>k\}|^{\alpha},
\end{equation}
where  $\{|w|>k\}=\{x\in\Omega:\, |w(x)|>k\}$.
Then $w\in L^\infty(\Omega)$ and
\begin{equation}\label{2.4}
 \norm{w}_{L^\infty(\Omega)}\le c_6 A^{\frac1\alpha} \norm{w}_{L^1(\Omega)}^{1-\frac1\alpha}.
\end{equation}
\end{lemma}
\noindent{\bf Proof.} We will follow the same idea of the proof  Lemma 5.2 in \cite{P}. For the convenience of the reader, we provide all the details of the proof.

By Cavalieri's principle, we have that
$$\norm{S_k(w)}_{L^1(\Omega)}=\int_0^\infty|\{|S_k(w)|>s\}|ds =\int_k^\infty|\{|w|>s\}|ds $$
Using (\ref{2.3}), it follows that
$$\int_k^\infty|\{|w|>s\}|ds\le A|\{|w|>k\}|^\alpha. $$
Let $H:[0,\infty)\to\R$ be the function defined for $t\ge0$ by
$$H(t)=\int_t^\infty|\{|w|>s\}|ds.$$
Since $s\to |\{|w|>s\}|$ is non-increasing, it is then continuous except for countable many points. Thus, for almost every $t\geq0$, we have that
$$-H'(t)=|\{|w|>t\}|\geq \left[\frac{H(t)}A\right]^{\frac1\alpha}.$$
Integrating this inequality, we conclude that if $\alpha>1$, then $H(k_0)=0$ for some $k_0\geq0$ such that
$$k_0\le c_7 A^{\frac1\alpha}H(0)^{1-\frac1\alpha}, $$
 that is, $\norm{S_k(w)}_{L^\infty(\Omega)}\le k_0$ and $H(0)=\norm{w}_{L^1(\Omega)}$. This  ends the proof. \hfill$\Box$\medskip

\begin{proposition}\label{pr 2.1}
Assume that $N\geq 3$,   $\mu>\frac34\mu_0$,   $f \in L^r(\Omega)$, $F  \in  L^r(\Omega;\R^N)$ with $r>N$ and $u_0\in \mathbb{H}^1_{0}(\Omega)$ is the unique solution of problem (\ref{eq 2.1}).

Then
 \begin{equation}\label{2.1}
   \norm{u_0}_{L^\infty(\Omega)}\le c_8\left(\norm{f }_{L^r(\Omega )} +\norm{ F }_{L^r(\Omega )} \right),
 \end{equation}
 where $c_8>0$ is independent of $f$ and $F$.
\end{proposition}
{\bf Proof.}      From (\ref{2.1-b5}) with $v=  S_k(u_0)\in \mathbb{H}^1_{0}(\Omega)\cap \mathbb{H}^1_{0}(\{|u_0|>k\})$, we have that
\begin{eqnarray}\label{cri identity-1}
 \int_{\{|u_0|>k\}} |\nabla u_0|^2 dx -2\tau_+\int_{\{|u_0|>k\}} \frac{x\cdot \nabla S_k^2(u_0)}{|x|^2} dx=\int_{\{|u_0|>k\}} \Big(f +  {\rm div}F \Big)    S_k(u_0)\, dx,
\end{eqnarray}
where we used the facts that $\nabla   S_k(u_0)=\nabla u_0$ a.e. in $\{|u_0|>k\}$ and $\nabla  S_k(u_0)=0$ a.e. in $\Omega\setminus \{|u_0|>k\}$.\smallskip

For  $\mu\geq 0$, it follows from (\ref{div 1}) that
\begin{eqnarray*}
 &&\int_{\{|u_0|>k\}} |\nabla u_0|^2 dx -2\tau_+\int_{\{|u_0|>k\}} \frac{x\cdot \nabla S_k^2(u_0)}{|x|^2} dx
 \\[1mm]&\geq &   \int_{\{|u_0|>k\} } |\nabla S_k(u_0)|^2 dx +(N-2)\tau_+  \int_{\{|u_0|>k\}} \frac{ S_k^2(u_0)}{|x|^2} dx
 \\[1mm]&\geq & \int_{\{|u_0|>k\}} |\nabla S_k(u_0) |^2 dx,
\end{eqnarray*}
where $\tau_+\geq0$ for $\mu\geq 0$.

For $0>\mu>-\frac34\frac{(N-2)^2}{4}$,  $\tau_+<0$ and  $1+\frac{4\tau_+}{N-2}>0$. It follows by    Hardy inequality and the fact that $ S_k(u_0)\in  \mathbb{H}^1_{0}(\{|u_0|>k\})$,  we have  that
$$\int_{\{|u_0|>k\} }\frac{S_k^2(u_0)}{|x|^2}  dx=\int_{\Omega}\frac{S_k^2(u_0)}{|x|^2}  dx\leq \frac4{(N-2)^2} \int_{\Omega} |\nabla S_k(u_0)|^2 dx$$
and
\begin{eqnarray*}
 &&\int_{\{|u_0|>k\}} |\nabla u_0|^2 dx -2\tau_+\int_{\{|u_0|>k\}} \frac{x\cdot \nabla S_k^2(u_0)}{|x|^2} dx
\\[1mm]&\geq&   \int_{\{|u_0|>k\} } |\nabla S_k(u_0)|^2 dx +(N-2)\tau_+ \int_{\{|u_0|>k\}} \frac{S_k^2(u_0)}{|x|^2}  dx
 \\[1mm]&=&(1+\frac{4\tau_+}{N-2})\int_{\{|u_0|>k\}} |\nabla S_k(u_0)|^2 dx.
\end{eqnarray*}

 The right hand side of (\ref{cri identity-1}) has the following estimate:
\begin{eqnarray*}
&& \Big|\int_{\{|u_0|>k\}} \Big(f  +  {\rm div}F\Big)     S_k(u_0)\, dx \Big|
\\[1mm]  &\leq &  (\int_{\{|u_0|>k\}}f ^2 \,dx )^{\frac12}  (\int_{\{|u_0|>k\}} S_k(u_0)^2dx )^{\frac12}  +\int_{\{|u_0|>k\}} |  F\cdot \nabla S_k(u_0)|\, dx
    \\[1mm]  & \leq &\Big(\norm{f }_{L^2(\{|u_0|>k\})}  + \norm{ F }_{L^2(\{|u_0|>k\})}\Big) \norm{\nabla S_k(u_0)}_{L^2(\{|u_0|>k\})}.
\end{eqnarray*}
It then follows from (\ref{cri identity-1}) that
\begin{equation}\label{2.2}
 \norm{\nabla S_k(u_0)}_{L^2(\{|u_0|>k\})}  \le  c_7 \Big(\norm{f }_{L^2(\{|u_0|>k\})} +\norm{F  }_{L^2(\{|u_0|>k\})}\Big).
\end{equation}
Therefore, it follows from  the H\"{o}lder inequality and  Sobolev inequality  that
\begin{eqnarray}
\norm{   S_k(u_0)  }_{L^1(\Omega)}&=&\norm{ S_k(u_0)}_{L^1(\{|u_0|>k\})}\nonumber
 \\ &\le &  \norm{S_k(u_0)}_{L^{2^*}(\Omega)}\, |\{|u_0|>k\} |^{\frac12+\frac1N}\nonumber \\
    &\le &\int_{\{|u_0|>k\} } |\nabla u_0|^2 dx\, |\{|u_0|>k\} |^{\frac12+\frac1N}   \label{3.1}  \\
    &\le &   c_9 \left(\norm{f }_{L^2(\{|u_0|>k\})} +\norm{F }_{L^2(\{|u_0|>k\})}\right)\, |\{|u_0|>k\}|^{\frac12+\frac1N}
   \nonumber \\&\le &  c_9 \left(\norm{f }_{L^r(\{|u_0|>k\})} +\norm{F  }_{L^r(\{|u_0|>k\})}\right)\, |\{|u_0|>k\}|^{1+\frac1N-\frac1r} ,\nonumber
\end{eqnarray}
that is,
\begin{equation}\label{2.02}
\norm{  S_k(u_0)}_{L^1(\Omega)}\le c_{10}\left(\norm{f }_{L^r(\{|u_0|>k\})} +\norm{F  }_{L^r(\{|u_0|>k\})}\right)\,|\{|u_0|>k\}|^{1+\frac1N-\frac1r},
\end{equation}
where $r>N$.\smallskip

  We apply Lemma \ref{lm 2.0} with  $\alpha=1+\frac1N-\frac1r>1$ for $r>N$, then it follows that
$$\norm{u_0 }_{L^\infty(\Omega)}\le c_{11} \left(\norm{f }_{L^r(\Omega)} +\norm{ F }_{L^r(\Omega)}\right), $$
which ends the proof. \hfill$\Box$

\begin{proposition}\label{pr 2.2}
Assume that $N\geq 3$,   $\mu>0$,   $f \in L^r(\Omega, , d\gg_\mu)$, $F  \in  L^r(\Omega, d\gg_\mu;\R^N)$ with $r>N+\tau_+$ and $u_0\in \mathbb{H}^1_{0}(\Omega)$ is the unique solution of problem (\ref{eq 2.1}).

Then
 \begin{equation}\label{2.1-1}
   \norm{u_0}_{L^\infty(\Omega)}\le c_{12}\left(\norm{f }_{L^r(\Omega, d\gg_\mu )} +\norm{ F }_{L^r(\Omega , d\gg_\mu)} \right),
 \end{equation}
 where $c_{12}>0$ is independent of $f$ and $F$.
\end{proposition}
{\bf Proof.} Let $D$ be a nonempty open domain, then  by H\"older inequality, we have that
\begin{eqnarray*}
 \norm{g }_{L^2(\CD)} ^2\leq   \norm{g }_{L^r(\CD, d\gg_\mu)}^2\, (\int_{\CD} \Gamma_\mu^{-\frac2r q_1} dx)^{\frac1{q_1}} \, |\CD|^{\frac1{q_2}},
\end{eqnarray*}
where our purpose is to find $r>2$ such that there exist $q_1,\, q_2>1$ satisfying that
\begin{equation}\label{3.1-1}
 \left\{\arraycolsep=1pt
\begin{array}{lll}
 \frac2r+\frac1{q_1}+\frac1{q_2}=1,\\[1.5mm]
 \phantom{ -\ \,\, }
 -\frac{2\tau_+}r q_1 >-N,\\[1.5mm]
 \phantom{ ----\ \, }
  \frac1{2q_2}>\frac12-\frac1N.
 \end{array}
 \right.
\end{equation}
Thanks to the fact $\tau_+>0$, when $r>N+\tau_+$, (\ref{3.1}) holds for some $q_1,\, q_2>1$.
Additionally, we note that for $r>N+\tau_+$
$$L^r(\CD,d\gg_\mu)\subset L^2(\CD).$$

Thus, from (\ref{3.1}), we have that
\begin{eqnarray*}
\norm{   S_k(u_0)  }_{L^1(\Omega)}&\leq &   c_{13} \left(\norm{f }_{L^2(\{|u_0|>k\})} +\norm{F }_{L^2(\{|u_0|>k\})}\right)\, |\{|u_0|>k\}|^{\frac12+\frac1N}
 \\[1mm]&\le &  c_{14} \left(\norm{f }_{L^r(\Omega, d\gg_\mu)} +\norm{F  }_{L^r(\Omega, d\gg_\mu)}\right)\, |\{|u_0|>k\}|^{\frac12+\frac1N+\frac1{q_2}},
\end{eqnarray*}
where $\frac12+\frac1N+\frac1{q_2}>1$. Then by Lemma \ref{lm 2.0} with  $\alpha=\frac12+\frac1N+\frac1{q_2}>1$ for $r>N$, then it follows that
$$\norm{u_0 }_{L^\infty(\Omega)}\le c_{15} \left(\norm{f }_{L^r(\Omega, d\gg_\mu)} +\norm{ F }_{L^r(\Omega, d\gg_\mu)}\right), $$
which ends the proof.\hfill$\Box$

\setcounter{equation}{0}
\section{ Proof of our main results}

From the form of $\nu=\gn\lfloor_{\Gw^*}+k\gd_0$ and linearity, we have to consider the regularity of the solutions of
(\ref{eq 1.1}) with $\nu=\gn\lfloor_{\Gw^*}$ and $\nu=k\gd_0$ respectively.

Let $v_0$ be a weak solution of
\begin{equation} \label{eq 3.1}
 \arraycolsep=1pt\left\{
\begin{array}{lll}
 \displaystyle    \mathcal{L}_\mu u = \delta_0 \quad\
   &{\rm in}\ \ {\Omega},\\[1.5mm]
 \phantom{   L_\mu \,   }
 \displaystyle  u= 0\quad \ &{\rm   on}\ \, \partial{\Omega},
 \end{array}\right.
 \end{equation}
 in the distributional sense that
\begin{equation}\label{weak sense 0}
\int_\Omega v_0 \mathcal{L}^*_\mu\xi \,d{\gg_\mu}(x)=c_\mu \xi(0)\quad \text{for all }\;\xi\in \BBX_\gm(\Gw),
\end{equation}
where $\mathcal{L}^*_\mu$ is given by (\ref{L}) and $c_{\mu}$ is the normalized constant.

\begin{proposition}\label{lm 3.1}
Let $N\geq 2$,  $\mu> \mu_0$  $p^*_\mu$ is given in (\ref{critical 1}) and $v_0$ be the very weak solution of   (\ref{eq 3.1}). 

 Then $(i)$     $v_0\Gamma_\mu \in W^{1,q}_0(\Omega)$ with $q\in[1, \frac{N}{N-1})$, 
\begin{equation}\label{5.2-0}
\norm{v_0\Gamma_\mu  }_{W^{1,q}(\Omega)} \le c_{16};
\end{equation}
$(ii)$     $v_0  \in W^{1,q}_0(\Omega,d\gg_\mu)$ with $q\in[1, p^*_\mu)$,
\begin{equation}\label{5.2-1}
\norm{v_0  }_{W^{1,q}(\Omega, d\gg_\mu)} \le c_{17}.
\end{equation}
\end{proposition}
{\bf Proof.} The existence of solution $v_0$  could see \cite[Theorem 1.2]{CQZ}.

Let $\eta_0:[0,+\infty)\to [0,\,1]$ be a decreasing $C^\infty$ function such that
\begin{equation}\label{eta}
 \eta_0=1\quad{\rm in}\quad [0,1]\qquad{\rm and}\qquad \eta_0=0\quad {\rm in}\quad[2,+\infty).
\end{equation}
Take $n_0\geq 1$ such that
$$\frac1{n_0}\sup\{r>0:\, B_r(0)\subset \Omega\}\le \frac12.$$
Denote $\eta_{n_0}(r)=\eta_0(n_0r)$ for $r\geq 0$,   $w_1=\Phi_\mu \eta_{n_0}$,
then Direct computation shows that
$w_1$ satisfies (\ref{5.2-0}) and (\ref{5.2-1}).

Note that 
\begin{equation}\label{6.1}
  \mathcal{L}_\mu w_1= -\nabla \eta_{n_0}\cdot \nabla  \Phi_\mu-\Phi_\mu \Delta\eta_{n_0}\quad
   {\rm in}\ \,  {\Omega}\setminus \{0\},
\end{equation}
where $-\nabla \eta_{n_0}\cdot \nabla  \Phi_\mu-\Phi_\mu \Delta\eta_{n_0}$ is smooth and has compact support in  $\overline{B_{\frac2{n_0}}(0)}\setminus B_{\frac1{n_0}}(0)$.

Let   $w_2$ be a solution of
\begin{equation}\label{6.2}
\arraycolsep=1pt\left\{
\begin{array}{lll}
 \displaystyle  \mathcal{L}_\mu w_2= -\nabla \eta_{n_0}\cdot \nabla  \Phi_\mu-\Phi_\mu \Delta\eta_{n_0}\qquad
   {\rm in}\quad  {\Omega}\setminus \{0\},\\[2mm]
 \phantom{  L_\mu \, }
 \displaystyle  w_2= 0\qquad  {\rm   on}\quad \partial{\Omega},\\[2mm]
 \phantom{   }
  \displaystyle \lim_{x\to0}w_i(x)\Phi_\mu^{-1}(x)= 0.
 \end{array}\right.
\end{equation}
then $w_2\in \BBH^1_0(\Omega)\cap C^2(\Omega\setminus\{0\})$ and 
$$\int_\Omega |\nabla w_2|^2dx<c_{19}.$$
 By comparison principle, we have that $|w_2|\leq c_{18}\Gamma_\mu$. 
 As a consequence, $w_2$ satisfies (\ref{5.2-0}) and (\ref{5.2-1}). Note that $v_0=w_1-w_2$ and the proof is complete.\hfill$\Box$

\medskip

\begin{lemma}\label{lm 3.2}
Assume that $N\geq 3$, $\mu>\frac{3}{4}\mu_0$ and $\nu\in \mathfrak{M}_+(\Omega^*; \Gamma_\mu)$.  Let $v_1$ be the unique weak solution of $\mathcal{L}_\mu   u =\nu \ \,  {\rm in}\ \,\Omega$ under the zero Dirichlet boundary condition.

Then    $v_1\Gamma_\mu \in W^{1,q}_0(\Omega)$ with $q\in[1,\frac{N}{N-1})$ such that
\begin{equation}\label{1.2-0}
\norm{v_1\Gamma_\mu  }_{W^{1,q}(\Omega)} \le c_{20} \norm{\nu}_{\overline{\mathfrak{M}}(\Omega; \Gamma_\mu)},
\end{equation}
where $c_{20}>0$ is independent of $\nu$.
\end{lemma}
\noindent{\bf Proof.}
For every $\xi\in C^{1.1}_0(\Omega)$, the $d\mu$-very  weak solution $u$ of (\ref{eq 1.1}) satisfies the inequality
\begin{equation}\label{ind 1}
\left|\int_\Omega u   \mathcal{L}^*_\mu(\xi) \, d\gamma_\mu \right|=\left|\int_\Omega \xi  \, d(\Gamma_\mu\nu) \right|\le \norm{\nu}_{\mathcal{M}(\Omega,\Gamma_\mu)} \norm{\xi  }_{L^\infty(\Omega)},\end{equation}
 where $d\gg_\mu=\Gamma_\mu dx$ and  $ d(\Gamma_\mu\nu)=\Gamma d\nu$.

For any $f\in C^1(\bar\Omega)$, let $\xi$ be the solution of (\ref{eq 2.1}) with $F=0$. For $q\in(1,\frac{N}{N-1})$, $q'=\frac{q}{q-1}>N$ and by Proposition \ref{pr 2.1}, it follows that
$$\norm{\xi }_{L^\infty(\Omega)}\le c_8\norm{f }_{L^{q'}(\Omega)}.$$
Therefore, we have that
$$\left|\int_\Omega u  \Gamma_\mu f\, dx \right| \le c_{8} \norm{\nu}_{\mathcal{M}(\Omega,\,\Gamma_\mu )} \norm{f }_{L^{q'}(\Omega)}. $$
Since $C^1(\bar \Omega)$ is dense in $L^{q'}(\Omega)$, then Riesz representation theorem implies that $u\Gamma_\mu \in L^q(\Omega)$ and
$$\norm{u\Gamma_\mu }_{L^q(\Omega)}\le c_{21}\norm{\nu}_{\mathcal{M}(\Omega,\,\Gamma)}.$$

For $F\in C^2(\bar\Omega;\R^N)$,  we may let $\xi$ be the solution of (\ref{eq 2.1}) with $f=0$.  By Proposition \ref{pr 2.1}, it follows that for $q\in(1,\frac{N}{N-1})$,
$$\norm{\xi }_{L^\infty(\Omega)}\le c_{8}\norm{F}_{L^{q'}(\Omega )}.$$
Henceforth, we have that
$$\left|\int_\Omega u\Gamma_\mu    {\rm div}F\, dx \right|\le c_{8} \norm{\nu}_{\mathcal{M}(\Omega,\,\Gamma_\mu )} \norm{F }_{L^{q'}(\Omega )}. $$
Consequently,  the functional
$$ F\in C^2(\bar\Omega;\R^N)\mapsto \int_\Omega u\Gamma_\mu \, {\rm div} F\, dx $$
admits a unique continuous linear extension in $L^{q'}(\Omega;\R^N)$. By the Riesz representation theorem, there exists a unique function $\tilde{\bf u}\in L^{q}(\Omega;\R^N)$
such that for any $F \in L^{q'}( \Omega;\R^N)$,
$$ \int_\Omega u \Gamma_\mu \,{\rm div}  F\, dx =\int_\Omega \tilde{\bf u}\cdot  F \,dx, $$
which implies that $\tilde{\bf u}=\nabla (u\Gamma_\mu)\in L^{q}(\Omega;\R^N). $
Thus, $u\Gamma_\mu  \in W^{1,q}_0(\Omega)$ and
$$\norm{u\Gamma_\mu  }_{W^{1,q}(\Omega)}\le c_{22} \norm{\nu}_{\mathcal{M}(\Omega,\,\Gamma_\mu)}$$
for $q\in (1,\frac{N}{N-1})$.    Theorem \ref{teo 1} with $q=1$  follows by H\"{o}lder inequality.\hfill$\Box$\medskip

\begin{lemma}\label{lm 3.3}
Assume that $N\geq 3$, $\mu>\frac34\mu_0$ $\mu\not=0$,  $\nu\in \mathfrak{M}_+(\Omega^*; \Gamma_\mu)$ and $p^*_\mu $ is given by (\ref{critical 1}).  Let $v_1$ be the unique weak solution of (\ref{eq 1.1})

Then    $v_1  \in W^{1,q}_0(\Omega, d\gamma_\mu)$ with $q\in[1, \min\{p^*_\mu,\frac{N}{N-1}\})$ such that
\begin{equation}\label{1.2-1}
\norm{v_1 }_{W^{1,q}(\Omega, d\gamma_\mu)} \le c_{23} \norm{\nu}_{ \mathfrak{M}(\Omega^*; \Gamma_\mu)},
\end{equation}
where $c_{23}>0$ is dependent of $\nu$.
\end{lemma}
\noindent{\bf Proof.} {\it The case of $\mu>0$.}
 For any $f\in C^1(\bar\Omega)$, let $\xi$ be the solution of (\ref{eq 2.1}) with $F=0$. For $q\in(1,\, p^*_\mu)$, $q'=\frac{q}{q-1}>N+\tau_+$ and by Proposition \ref{pr 2.2}, it follows that
$$\norm{\xi }_{L^\infty(\Omega)}\le c_8\norm{f }_{L^{q'}(\Omega,d\gamma_\mu)}.$$
Therefore, (\ref{ind 1}) implies that
$$\left|\int_\Omega u  \Gamma_\mu f\, dx \right| \le c_8 \norm{\nu}_{\mathcal{M}(\Omega,\,\Gamma_\mu )} \norm{f }_{L^{q'}(\Omega, d\gg_\mu)}. $$
Since $C^1(\bar \Omega)$ is dense in $L^{q'}(\Omega)$, then Riesz representation theorem implies that $u  \in L^q(\Omega, d\gamma_\mu)$ and
$$\norm{u }_{L^q(\Omega, d\gamma_\mu)}\le c_8\norm{\nu}_{\mathcal{M}(\Omega,\,\Gamma_\mu)}.$$

For $F\in C^2(\bar\Omega;\R^N)$,  we may let $\xi$ be the solution of (\ref{eq 2.1}) with $f=0$.  By Proposition \ref{pr 2.2}, it follows that for $q\in(1,p^*_\mu)$,
$$\norm{\xi }_{L^\infty(\Omega)}\le c_{17}\norm{F}_{L^{q'}(\Omega, d\gamma_\mu )}.$$
Henceforth, we have that
$$\left|\int_\Omega u\Gamma_\mu    {\rm div}F\, dx \right|\le c_{18} \norm{\nu}_{\mathcal{M}(\Omega,\,\Gamma_\mu )} \norm{F }_{L^{q'}(\Omega , d\gamma_\mu)}. $$
Consequently,  the functional
$$ F\in C^2(\bar\Omega, d\gamma_\mu;\R^N)\mapsto \int_\Omega u\Gamma_\mu \, {\rm div} F\, dx $$
admits a unique continuous linear extension in $L^{q'}(\Omega, d\gamma_\mu;\R^N)$. By the Riesz representation theorem, there exists a unique function $\tilde{\bf u}\in L^{q}(\Omega,d\gg_\mu;\R^N)$
such that for any $F \in L^{q'}(\bar\Omega, d\gg_\mu;\R^N)$,
$$ \int_\Omega u \Gamma_\mu \,{\rm div}  F\, dx =\int_\Omega \tilde{\bf u}\cdot  F \,dx=\int_\Omega (\gamma_\mu^{-\frac1{q'}}\tilde{\bf u})\cdot  (\gamma_\mu^{\frac1{q'}} F) \,dx, $$
which implies that $\Gamma_\mu^{-\frac1{q'}}\tilde{\bf u}=\Gamma_\mu^{-\frac1{q'}} \nabla (u\Gamma_\mu)=\Gamma_\mu^{\frac1{q}}\nabla u+\Gamma_\mu^{\frac1{q}}u\, \frac{x}{|x|^2}\in L^{q}(\Omega;\R^N)$ and
$$\norm{(\Gamma_\mu^{\frac1{q}}\nabla u+\Gamma_\mu^{\frac1{q}}u\, \frac{x}{|x|^2} ) }_{L^q(\Omega)}\le c_{23} \norm{\nu}_{\mathcal{M}(\Omega,\,\Gamma_\mu)}.$$

{\bf  Claim 1: } for $1\leq q<p^*_\mu$, we have that
$$\|\Gamma_\mu^{\frac1{q}}u\, \frac{x}{|x|^2} \|_{L^q(\Omega)}\leq c_{24}\norm{\nu}_{ \mathfrak{M}(\Omega^*; \Gamma_\mu)}.$$

 Thus,
 $\Gamma_\mu^{\frac1{q}}\nabla u=[(\Gamma_\mu^{\frac1{q}}\nabla u+\Gamma_\mu^{\frac1{q}}u\, \frac{x}{|x|^2} )-\Gamma_\mu^{\frac1{q}}u\, \frac{x}{|x|^2}]\in L^q(\Omega;\R^N)$,
 i.e. $|\nabla u|\in L^q(\Omega,\,d\gamma_\mu)$.\smallskip

We end the proof by showing {\bf  Claim 1}.  In fact, by Lemma \ref{lm 3.2}, we have that $u\Gamma_\mu\in W^{1,p}_0(\Omega)$ for any $p\in [1,\, \frac{N}{N-1})$, then by Sobolev Embedding theorem, for any $\gs\in[1, \frac{N}{N-2})$, there holds that
$$\|u\Gamma_\mu\|_{L^\gs(\Omega)}\leq c_{25}\norm{\nu}_{ \mathfrak{M}(\Omega^*; \Gamma_\mu)}.$$

We observe that for $1\leq q<p^*_\mu $,
\begin{eqnarray*}
\int_\Omega \Big|\Gamma_\mu^{\frac1{q}}u\, \frac{x}{|x|^2}\Big|^q dx&= & \int_\Omega |x|^{\tau_+(1-q)-q} | |x|^{\tau_+} u|^q  dx
 \\[1mm] &\leq & (\int_\Omega | |x|^{\tau_+} u|^{q_1}  dx)^{\frac{q}{q_1}} (\int_{\Omega} |x|^{[\tau_+(1-q)-q]q_2}dx)^{\frac1{q_2}},
\end{eqnarray*}
where $q_2=\frac{q_1}{q_1-q}$ and we require that
\begin{equation}\label{12-1}
q_1<\frac{N}{N-2}
\end{equation}
and
\begin{equation}\label{12-1a1}
[\tau_+(1-q)-q]q_2>-N.
\end{equation}
(\ref{12-1a1}) could be equivalent to
$ Nq< [N-q-\tau_+(q-1)]q_1$, then by (\ref{12-1}), we have that
$$ (N-2)q<N-q-\tau_+(q-1),$$
 which holds for any $1\leq q<p^*_\mu $.

{\it For the case $\mu\in(\frac34\mu_0,0)$,}  we note that
\begin{eqnarray*}\int_{\Omega}|\nabla u|^pd\gg_\mu &= &  \int_{\Omega}\Gamma_\mu |\nabla (\Gamma_\mu^{-1}(\Gamma_\mu u))|^pdx
\\&\leq& 2^p\Big(|\tau_+|^p \int_{\Omega} |x|^{(1-p)\tau_+-p} |\Gamma_\mu u|^p dx+\int_{\Omega}|x|^{(1-p)\tau_+} |\nabla (\Gamma_\mu u)|^p \Big)
\\&\leq& c_{26} \Big(\int_{\Omega}|\Gamma_\mu u|^p dx+\int_{\Omega}|\nabla(\Gamma_\mu u)|^p dx\Big)
\\&\leq& c_{27}\norm{\nu}_{ \mathfrak{M}(\Omega^*; \Gamma_\mu)} ^p,
 \end{eqnarray*}
where we used  $p\geq 1$ and $\tau_+<0$ for $\tau<0$.  We complete the proof. \hfill$\Box$\medskip\smallskip

 \noindent{\bf Proof of Theorem \ref{teo 1}.} Part $(i)$  follows directly by Proposition \ref{lm 3.1} and Lemma \ref{lm 3.2};
 Part $(ii)$ does Proposition \ref{lm 3.1} and Lemma \ref{lm 3.3}. \hfill$\Box$ \medskip

  \noindent{\bf Proof of Corollary \ref{cr 1}.} Note that from definition \ref{weak definition},  any weak solution $u$ of (\ref{eq 1.1-absorption})
   has the property that $g(u)\in L^1(\Omega, d\gg_\mu)$. In fact,
 for any $w\in L^1(\Omega, d\gg_\mu)$,
we have that
 $$\int_{\Omega^*} |w| d \gg_\mu= \int_{\Omega} |w| d\gg_\mu<+\infty,$$
 then
  $L^1(\Omega, d\gg_\mu)\subset
  \mathfrak{M}(\Omega^*; \Gamma_\mu)$ and then
    $\nu-g(u)\in \overline{\mathfrak{M}}(\Omega; \Gamma_\mu)$. Apply Theorem \ref{teo 1}, we obtain the global regularity.\hfill$\Box$

    \section{More global regularity: the Marcinkiewicz estimates}

In this section, we append  a bit discussion of global regularity for non-homogeneous 
Hardy problem (\ref{eq 1.1}). To this end, we recall the definition and basic properties of the Marcinkiewicz
spaces.

\begin{definition}\label{definition 2.1}
Let $\Omega\subset \R^N$ be a domain and $\nu$ be a positive
Borel measure in $\Omega$. For $\kappa>1$,
$\kappa'=\kappa/(\kappa-1)$ and $u\in L^1_{loc}(\Omega,d\nu)$, we
set
$$
 \arraycolsep=1pt
\begin{array}{lll}
\|u\|_{M^\kappa(\Omega,d\nu)}=\inf\{c\in[0,\infty]:\int_E|u|d\nu\le
c\left(\int_Ed\nu\right)^{\frac1{\kappa'}},\ \forall E\subset\Omega\
\rm{Borel\ set}\}
\end{array}
$$
and
\begin{equation}\label{a.spa M}
M^\kappa(\Omega,d\nu)=\{u\in
L_{loc}^1(\Omega,d\nu):\|u\|_{M^\kappa(\Omega,\,d\nu)}<+\infty\}.
\end{equation}
\end{definition}
$M^\kappa(\Omega,d\nu)$ is called the Marcinkiewicz space with
exponent $\kappa$ or weak $L^\kappa$ space and
$\|.\|_{M^\kappa(\Omega,d\nu)}$ is a quasi-norm. The following
property holds.

\begin{proposition}\label{a.pr 1} \cite{BBC,CC1}
Assume that $1\le q< \kappa<+\infty$ and $u\in L^1_{loc}(\Omega,d\nu)$.
Then there exists  $C(q,\kappa)>0$ such that
$$\int_E |u|^q d\nu\le C(q,\kappa)\|u\|_{M^\kappa(\Omega,d\nu)}\left(\int_E d\nu\right)^{1-q/\kappa}$$
for any Borel set $E$ of $\Omega$.
\end{proposition}


\begin{proposition}\label{a.general}
Let $\mu>0$, $\Omega\subset \R^N\ (N\geq2)$ be a bounded $C^2$ domain
and $\mathbb{G}_\mu[\nu](\cdot)=\int_{\Omega}G_\mu(\cdot,y)d\sigma(y)$ for  $\sigma\in\mathfrak{M}(\Omega,\Gamma_\mu)$. Then there exists $c_{28}>0$ such that
\begin{equation}\label{a.annex 0}
\|\mathbb{G}_\mu[\sigma]  \|_{M^{q}(\Omega,\,d\gamma_\mu )}\le c_{28}\|\sigma\|_{\mathfrak{M}(\Omega,\, \Gamma_\mu)},
\end{equation}
where $1<q\le \frac{N+\tau_+(\mu)}{N-2+\tau_+(\mu)}$.
\end{proposition}

\noindent{\bf Proof.}  For $\lambda>0$ and $y\in\Omega$, we denote
$$A_\lambda(y)=\{x\in\Omega\setminus\{y\}: G_\mu(x,y)>\lambda\}\ \ {\rm
{and}}\quad m_\lambda(y)=\int_{A_\lambda(y)} d\gamma_\mu.$$
From \cite[Lemma 4.1]{CQZ} if follows  that for any $(x,y)\in
(\Omega\setminus\{0\})\times(\Omega\setminus\{0\})$, $x\neq y$,
\begin{eqnarray*}
G_\mu(x,y) &\le& c_{29}  \frac{\Gamma_\mu(y)}{|x-y|^{N-2+\tau_+(\mu)}}.
\end{eqnarray*}

Observe that
\begin{eqnarray*}
A_\lambda(y)\subset \left\{x\in\Omega\setminus\{y\}: \frac{c_{29}\Gamma_\mu(y)}{|x-y|^{N-2+\tau_+(\mu)}} >\lambda\right\}=  D_\lambda(y),
\end{eqnarray*}
where
$D_\lambda(y)=\left\{x\in\Omega:|x-y|<(\frac{c_{29}\Gamma_\mu(y)}{\lambda})^{\frac1{N-2+\tau_+(\mu)}}\right\}$.
Then we have that
\begin{equation}\label{a.annex 1xhw}
m_\lambda(y)\le\int_{D_\lambda(y)}  d\gamma_\mu\le\int_{D_\lambda(0)}  d\gamma_\mu 
= c_{29}(\frac{c_{29}\Gamma_\mu(y)}{\lambda})^{\frac {N+\tau_+(\mu)}{N-2+\tau_+(\mu)}}.
\end{equation}
For any Borel set $E$ of $\Omega$, we have
\begin{eqnarray*}
\int_E G_\mu(x,y)\Gamma_\mu(x)dx\le
\int_{A_\lambda(y)}G_\mu(x,y)\Gamma_\mu(x)dx+\lambda\int_E d\gamma_\mu.
\end{eqnarray*}
Using Fubini's theorem, integration by parts formula and estimate (\ref{a.annex 1xhw}), we obtain
\begin{eqnarray*}
\int_{A_\lambda(y)}G_\mu(x,y) d\gamma_\mu &=&-\int_\lambda^\infty sdm_s(y)
ds
\\&=&\lambda m_\lambda(y)+ \int_\lambda^\infty m_s(y)ds
\\&\le& c_{30}\Gamma_\mu^{\frac{N+\tau_+(\mu)}{N-2+\tau_+(\mu)}}(y)\lambda^{1-\frac{N+\tau_+(\mu)}{N-2+\tau_+(\mu)}}.
\end{eqnarray*}
Thus,
\begin{eqnarray*}
\int_E G_\mu(x,y) d\gamma_\mu\le
c_{30}\Gamma_\mu^{\frac{N+\tau_+(\mu)}{N-2+\tau_+(\mu)}}(y)\lambda^{1-\frac{N+\tau_+(\mu)}{N-2+\tau_+(\mu)}} +\lambda \int_E d\gamma_\mu.
\end{eqnarray*}
By choosing $\lambda=\Gamma_\mu(y) (\int_E d\gamma_\mu)^{-\frac{N-2+\tau_+(\mu)}{N+\tau_+(\mu)}}$ and $c_{31}=c_{30}+1$, we have
\begin{eqnarray*}
\int_E G_\mu(x,y)d\gamma_\mu \le c_{31} \Gamma_\mu(y)
(\int_E  d\gamma_\mu)^{\frac{2+\tau_+(\mu)}{N+\tau_+(\mu)}}.
\end{eqnarray*}
Therefore,
\begin{eqnarray*}
 \int_E
\mathbb{G}_\mu[|\sigma|](x)\Gamma_\mu(x)dx&=&\int_\Omega\int_E
G_\mu(x,y)\Gamma_\mu(x)dx d|\sigma(y)|
\\&\le &c_{31}\int_\Omega  \Gamma_\mu(y) d|\sigma (y)|\left(\int_E d\gamma_\mu\right)^{\frac{2+\tau_+(\mu)}{N+\tau_+(\mu)}}
\\&\le& c_{31}\|\sigma\|_{\mathfrak{M}(\Omega,\Gamma_\mu)} \left(\int_E d\gamma_\mu\right)^{\frac{2+\tau_+(\mu)}{N+\tau_+(\mu)}}.
\end{eqnarray*}
 As a consequence,
\begin{eqnarray*}
\|\mathbb{G}_\mu[\sigma]  \|_{M^{\frac{N+\tau_+(\mu)}{N-2+\tau_+(\mu)}}(\Omega,d\gamma_\mu)}\le
c_{31}\|\sigma\|_{\mathfrak{M}(\Omega,\Gamma_\mu)}.
\end{eqnarray*}
We complete the proof. \hfill$\Box$\smallskip

We remark that the biggest defect of this  Marcinkiewicz estimate in Proposition \ref{a.general} is the restriction
that $\mu>0$. For $N\geq 3$ and $\mu_0\leq \mu<0$, a much weaker type of Marcinkiewicz estimate could 
be derived due to the Green kernel's estimates, see \cite[Lemma 4.1]{CQZ}.

 \bigskip  \bigskip

 \noindent{\bf Acknowledgements:} H. Chen is supported by NSF of China, No:  11661045 and Key R$\&$D plan of Jiangxi Province, 20181ACE50029.

\end{document}